\documentclass[EJP]{ejpecp} 

\AUTHORS{%
  Xiequan Fan\footnote{Regularity Team, Inria and MAS Laboratory, Ecole Centrale Paris - Grande Voie des Vignes, 92295 Ch\^{a}tenay-Malabry, France. Corresponding author. E-mail: fanxiequan@hotmail.com}\ \
    \and
  Ion Grama\footnote{Universit\'{e} de Bretagne-Sud,  UMR 6205, LMBA, 56000 Vannes, France. E-mail: ion.grama@univ-ubs.fr, quansheng.liu@univ-ubs.fr}
  \and
  Quansheng Liu\footnotemark[2]}

\SHORTTITLE{Exponential inequalities for martingales with applications}

\TITLE{Exponential inequalities for martingales with~applications} 


\AMSSUBJ{60E15; 60G42} 

\SUBMITTED{May 6, 2014} 
\ACCEPTED{December 27, 2014} 

\VOLUME{20}
\YEAR{2015}
\PAPERNUM{1}
\DOI{v20-3496}

\ABSTRACT{The paper is devoted to establishing some general exponential
  inequalities for supermartingales. The inequalities improve or generalize
  many exponential inequalities of Bennett, Freedman, de la Pe\~{n}a, Pinelis
  and van de Geer. Moreover, our concentration inequalities also improve some
  known inequalities for sums of independent random variables. Applications
  associated with linear regressions, autoregressive processes and branching
  processes are provided. In particular, an interesting application of
  \mbox{de la Pe\~{n}a's} inequality to self-normalized deviations is also
  provided. }

\begin{document}



\section{Introduction}
Assume that we are given a sequence of real-valued supermartingale differences $(\xi _i,\mathcal{F}_i)_{i=0,...,n}$
 defined on some probability space $(\Omega ,\mathcal{F},\mathbf{P})$, where $\xi _0=0 $ and $\{\emptyset, \Omega\}=%
\mathcal{F}_0\subseteq ...\subseteq \mathcal{F}_n\subseteq \mathcal{F}$ are
increasing $\sigma$-fields. So we have $\mathbf{E}(\xi_{i}|\mathcal{F}_{i-1})\leq 0, \  i=1,...,n, $ by definition. Set
\begin{equation}  \label{matingal}
 S_k=\sum_{i=1}^k\xi _i,\quad k=1,...,n.
\end{equation}
Then $S=(S_k ,\mathcal{F}_k)_{k=1,...,n}$ is a supermartingale.
Let $\left\langle S\right\rangle $ and $[S]$ be respectively the quadratic characteristic and the squared variation of the
supermartingale $S:$
\begin{equation}
 \left\langle S\right\rangle
_k=\sum_{i=1}^k\mathbf{E}(\xi _i^2|\mathcal{F} _{i-1})\ \ \ \ \ \textrm{and} \ \ \ \ \ \ \  [S]_k=\sum_{i=1}^{k}\xi_i^2.
\end{equation}
The following exponential inequality for supermartingales can be found in Freedman \cite{FR75}. \\
\textbf{Theorem A.\ }
Suppose  $\xi_{i} \leq \epsilon$ for a positive constant $\epsilon$.   Then, for all $x, v > 0,$
\begin{eqnarray}
 \mathbf{P}\Big(S_k \geq x\ \mbox{and}\ \left\langle S\right\rangle_k\leq v^2\ \mbox{for some}\ k \Big)
 &\leq& B_{2}(x, \epsilon, v) \label{freedmanxf2} \\
 &:=& \exp \left\{-\frac{x^2}{ 2(v^2+ x \epsilon )} \right\}. \nonumber
\end{eqnarray}
After Freedman's seminal work, many interesting exponential inequalities for martingales have been established.
For continuous-time martingales with bounded jumps, \mbox{Freedman's} inequality (\ref{freedmanxf2}) has been established by Shorack and Wellner \cite{SW86}.
By imposing certain moment conditions, \mbox{van de Geer} \cite{V95} relaxed the condition of Shorack and Wellner and generalized  inequality (\ref{freedmanxf2}) for martingales with non-bounded jumps.  Under the following conditional Bernstein condition: for a positive constant $\epsilon$,
\begin{eqnarray}
\mathbf{E} (|\xi_{i}|^l|\mathcal{F} _{i-1})  \leq \frac{1}{2}\, l !\, \epsilon^{l-2}\mathbf{E} (\xi _i^2|\mathcal{F} _{i-1}),\ \ \  \mbox{for all} \ \ l\geq 2, \label{Bernstein}
\end{eqnarray}
de la Pe\~{n}a \cite{D99} have obtained the following Bernstein type inequality for martingales, for all $x, v > 0,$
\begin{eqnarray}
\mathbf{P}(S_k \geq x\ \mbox{and}\ \langle S\rangle_{k}\leq v^2\ \mbox{for some}\ k )  &\leq& B_{1}(x, \epsilon, v) \label{djknss}\\
&:=& \exp \bigg\{-\frac{x^2}{v^2(1+\sqrt{1+2 x \epsilon/v^2})+x \epsilon  } \bigg\}\nonumber \\
 &\leq& B_{2}(x, \epsilon, v).  \label{djkss}
\end{eqnarray}
Inequality (\ref{djkss}) has also been  obtained by \mbox{van de Geer} \cite{V95}. In  particular, when $(\xi_i)_{i=1,...,n}$ are independent, the inequalities (\ref{djknss}) and (\ref{djkss}) reduce, respectively, to the inequalities of Bennett \cite{B62}  and Bernstein  \cite{B27}.  Many other generalizations of Freedman's inequality can be found in  Haeusler  \cite{H84}, Pinelis \cite{P94a}, Dzhaparidze  and van Zanten \cite{Dz01}, Delyon \cite{De09}  and  Khan \cite{K09}.

Following the work of Freedman \cite{FR75}, Shorack and Wellner \cite{SW86},
van de Geer \cite{V95} and de la Pe\~{n}a \cite{D99}, we develop some new
methods, based on changes of probability measure, for establishing some
general  exponential inequalities for supermartingales. The
methods are user-friendly and efficient.

In  Theorem \ref{th5}, we obtain two  exponential inequalities for supermartingales under a very general condition.
Assume that $$ \mathbf{E}(\exp\left\{\lambda \xi_i - g(\lambda)\xi_i^2     \right\}|\mathcal{F}_{i-1}) \leq  1+ f(\lambda)V_{i-1}$$ for some  $\lambda \in (0, \infty),$ for two non-negative functions $ f(\lambda)$ and $g(\lambda)$,  and for some non-negative and $\mathcal{F}_{i-1}$-measureable random variables $V_{i-1}$.  Then, for all $ x, v, \omega > 0$,
\begin{eqnarray}
&& \mathbf{P}\Big( S_{k} \geq x,\ [S]_{k} \leq v^2 \ \mbox{and}\ \sum_{i=1}^{k} V_{i-1} \leq w  \ \mbox{for some}\ k \in [1, n]  \Big) \nonumber\\
&\leq&  \exp \left\{ -\lambda x + g(\lambda)v^2+n\log\left( 1+ \frac{f(\lambda)}{n} w \right)\right\}  \label{gnhs1}\\
&\leq& \exp \Big\{ -\lambda x + g(\lambda)v^2 +f(\lambda)  w   \Big\} . \label{gnhs2}
\end{eqnarray}

If $\xi_i \geq -\epsilon$ for a positive constant $\epsilon$, then our result (\ref{gnhs2}) implies that, for all $x, v >0$,
\begin{eqnarray}
 \mathbf{P}\left( S_k  \geq x\ \mbox{and}\ [S]_k\leq v^2 \ \mbox{for some}\ k   \right)
 \leq  B_2\left(x, \epsilon, v\right).\label{ghssfyl2}
\end{eqnarray}
This  inequality is similar
to the one of Freedman (\ref{freedmanxf2}).
To highlight the differences between (\ref{freedmanxf2})  and  (\ref{ghssfyl2}), notice that the conditions $\xi_i \leq \epsilon$ and conditional variance $\langle S\rangle_k$ in Freedman's inequality (\ref{freedmanxf2}) are respectively replaced by the condition $\xi_i\geq -\epsilon$  and  squared variation $[S]_k$ in our inequality (\ref{ghssfyl2}).
Moreover, inequality (\ref{ghssfyl2}) completes Freedman's inequality (\ref{freedmanxf2}) by giving an estimation of deviation probabilities on the left side:
if the martingale differences $ (\xi _i,\mathcal{F}_i)_{i=1,..., n }$ satisfy $\xi_i \leq \epsilon$ for all $i$, then, for all $x, v >0$,
\begin{eqnarray}
 \mathbf{P}\left( S_k  \leq -x\ \mbox{and}\ [S]_k\leq v^2 \ \mbox{for some}\ k   \right)  \leq  B_2\left(x, \epsilon, v\right).
\end{eqnarray}

If the martingale differences verifies canonical assumption (which means $g(\lambda)=\lambda^2/2$ and $f(\lambda)=0$), then (\ref{gnhs2}) implies the following de la Pe\~{n}a inequality \cite{D99}, for all $x, v>0$,
\begin{eqnarray}\label{dlpieq}
 \mathbf{P}\Big(   S_{k} \geq x\ \mbox{and}\  [S]_k \leq v^2\ \mbox{for some}\ k  \Big)  \leq  \exp \left\{ - \frac{x^2}{2\,v^2} \right\}.
\end{eqnarray}
Moreover, we find that  (\ref{dlpieq}) implies the following self-normalized deviation result associated with  independent and symmetric random variables, for all $x>0$,
\begin{eqnarray}
 \mathbf{P}\bigg(  \max_{1\leq k \leq n} \frac{ S_{k}}{\sqrt{[S]_n} } \geq x   \bigg)  \leq  \exp \left\{ - \frac{x^2}{2} \right\}.
\end{eqnarray}

If $\mathbf{E} |\xi _i|^3 < \infty$, then (\ref{gnhs2}) implies  the following Bernstein type inequality, for all $x, v, w>0$,
\begin{eqnarray}
 \mathbf{P}\Big( S_k  \geq x,\ [S]_k \leq v^2\ \mbox{and}\  \Upsilon(S_k)\leq w  \ \mbox{for some}\ k \Big) &\leq&  B_1\left(x, \frac{w }{3v^2}, v\right) \label{fineq2}\\
&\leq& B_2\left(x, \frac{w }{3v^2}, v\right), \label{fineq3}
\end{eqnarray}
where $$ \Upsilon(S_k) =\sum_{i=1}^k \mathbf{E}(|\xi_i|^3  |\mathcal{F}_{i-1});$$
see Corollary \ref{co9}. Compared to the inequalities (\ref{djknss}) and  (\ref{djkss}), the advantage of the last two inequalities (\ref{fineq2}) and (\ref{fineq3}) is that we do not assume the existence of moments of all orders.

Assume that $\mathbf{E}(e^{\lambda \xi_{i}}|\mathcal{F}_{i-1}) \leq  1+ f(\lambda)\mathbf{E}(\xi_i^2 |\mathcal{F}_{i-1})$ for some   $\lambda \in (0, \infty)$ and a positive function $f(\lambda).$  Then Theorem \ref{th5} implies that, for all $ x, v > 0$,
\begin{eqnarray}
 && \mathbf{P}\left( S_k \geq x\ \mbox{and}\ \langle S\rangle_{k}\leq v^2\ \mbox{for some}\ k \in [1, n] \right) \nonumber\\
&\leq&  \exp \left\{ -\lambda x +n\log\left( 1+ \frac{f(\lambda)}{n} v^2\right)\right\} \ \ \   \label{ff07} \\
&\leq& \exp \Big\{ -\lambda x +f(\lambda)  v^2 \Big\} \label{ff08}.
\end{eqnarray}
In particular,  if $(\xi _i,\mathcal{F}_i)_{i=1,...,n}$ satisfies condition (\ref{Bernstein}), then it holds
$$ \mathbf{E}(e^{\overline{\lambda} \xi_{i}}|\mathcal{F}_{i-1}) \leq  1+ f(\overline{\lambda})\mathbf{E}(\xi_i^2 |\mathcal{F}_{i-1}),$$ where $$ \overline{\lambda} = \frac{2x/v^2}{2x\epsilon/v^2+1+\sqrt{1+2x\epsilon/v^2}} \ \ \ \textrm{and}  \ \ \ f(\lambda)=\frac{\lambda^2   }{2 (1-\lambda\epsilon)}.$$ Inequality (\ref{ff08}) reduces to \mbox{de la Pe\~{n}a's} inequality (\ref{djknss}) with $\lambda = \overline{\lambda}$. Hence,  our bound (\ref{ff07}) with $\lambda = \overline{\lambda}$ improves \mbox{de la Pe\~{n}a's} inequality (\ref{djknss}). In the i.i.d.\ case, bound (\ref{ff07}) significantly improves the large deviation bound (\ref{djknss}) on large deviation  tail probabilities $\mathbf{P}(S_n \geq nx)$ by adding a  factor with  exponentially decay rate $\exp\{ -n c_{x}\},$ where $c_{x}>0$ does not depend on $n$. In the applications for linear regression models, we find that such type refinements are useful; see Theorem \ref{thlin}.

In Theorem \ref{th3}, we consider the case that supermartingale has sub-Gaussian differences. Assume that $ \mathbf{E}(e^{\lambda \xi_{i}}|\mathcal{F}_{i-1}) \leq  \exp \{ f(\lambda) V_{i-1}  \}$ for some $\lambda \in (0, \infty),$ for  a positive function $ f(\lambda)$ and for some $\mathcal{F}_{i-1}$-measurable random variables $V_{i-1}$. Then, for all $ x, v > 0$,
\begin{eqnarray}
 \mathbf{P}\Big( S_{k} \geq x\ \mbox{and}\ \sum_{i=1}^{k} V_{i-1} \leq v^2\ \mbox{for some}\ k   \Big)
  \leq  \exp \Big\{ -\lambda x +f(\lambda)  v^2 \Big\}. \label{fsat9}
\end{eqnarray}
In particular, when the function $f(\lambda)=\lambda^2/2$ for all $\lambda>0$ and  $(V_{i})_{i=1,..,n}$ are constants,
inequality (\ref{fsat9}) reduces to Fuk's inequality with $\lambda= \lambda(x):=x/\sum_{i=1}^n V_i^2$ (cf.\ Theorem 4 of \cite{F73}). Thus (\ref{fsat9}) is
a generalization of Fuk's inequality \cite{F73} for supermartingales. If $V_{i-1}=\mathbf{E}(\xi_i^2|\mathcal{F}_{i-1})$ is  the conditional variance, inequality (\ref{fsat9}) reduces to Theorem 4.2 of Khan \cite{K09}.
Inequality (\ref{fsat9}) implies the following result, where $V_{i-1}$ is not  the conditional variance.
If $\xi_{i} \leq U_{i-1}$ for some $\mathcal{F}_{i-1}$-measurable random variables $U_{i-1}$,
then, for all $x, v >0$,
\begin{eqnarray}\label{fint3}
 \mathbf{P}\bigg(  S_k \geq x\  \mbox{and}\  \sum_{i=1}^k C_{i-1}^2 \leq v^2\ \mbox{for some}\ k    \bigg)  \leq   \exp\left\{- \frac{ x^2}{  2\, v^2 }   \right\},
\end{eqnarray}
where
\begin{eqnarray}  \label{f1}
C_{i-1}^2 = \left\{ \begin{array}{ll}
\mathbf{E}(\xi_{i}^2|\mathcal{F}_{i-1}) , & \textrm{\ \ \ \ \ if $\mathbf{E}(\xi_{i}^2|\mathcal{F}_{i-1}) \geq U^2_{i-1}$ }, \\
\displaystyle\frac{1}{4}\left( U_{i-1}  + \frac{\mathbf{E}(\xi_{i}^2|\mathcal{F}_{i-1}) }{ U_{i-1} }\right)^2, & \textrm{\ \ \ \ \ otherwise}.
\end{array} \right.
\end{eqnarray}
Then we show that (\ref{fint3}) implies a generalization of Azuma-Hoeffding's
inequality for martingales due to van de Geer \cite{V02}. Moreover, we also
show that (\ref{fint3}) significantly improves some recent inequalities of
\mbox{Bentkus} \cite{Be03} and Pinelis \cite{P06,P06b} by adding an
exponential decay factor in the case of $||\sum_{i=1}^n C_{i-1}^2||_{\infty}$ $ <
\sum_{i=1}^n ||C_{i-1}^2||_{\infty}$; see (\ref{mspi}) and \emph{Example} 1
for details. We find that such improvements are important in the applications
for linear regression models and autoregressive processes; see Remarks
\ref{fsfdf} and \ref{fsfdf1}.

The paper is organized as follows. We present our theoretical results in Section \ref{sec2}, give the applications of our results in Section \ref{sec2.5} and devote to the proofs of our results in Sections  \ref{sec6} - \ref{endsec}. The proofs of the theorems and their corollaries are in the same sections.

\section{Main results}\label{sec2}
Our first result is given under a very general condition.
\begin{theorem}\label{th5} Assume that $V_{i-1}, i \in [1, n],$ are non-negative and $\mathcal{F}_{i-1}$-measureable random variables. Suppose that
\begin{equation} \label{condimenme1}
\mathbf{E}(\exp\left\{\lambda \xi_i - g(\lambda)\xi_i^2     \right\}|\mathcal{F}_{i-1}) \leq  1+ f(\lambda)V_{i-1}
\end{equation}
for some $\lambda\in (0, \infty)$, for two non-negative functions $f(\lambda)$ and  $g(\lambda)$, and for all $i \in [1, n]$.  Then, for all $ x, v, \omega > 0$,
\begin{eqnarray}
&& \mathbf{P}\bigg( S_{k} \geq x,\ [S]_{k} \leq v^2 \ \mbox{and}\ \sum_{i=1}^{k} V_{i-1} \leq w \ \mbox{for some}\ k \in [1, n]  \bigg) \nonumber\\
&\leq&  \exp \left\{ -\lambda x + g(\lambda)v^2+n\log\left( 1+ \frac{f(\lambda)}{n} w \right)\right\}  \label{fhgsa1} \\
&\leq& \exp \Big\{ -\lambda x + g(\lambda)v^2 +f(\lambda)  w  \Big\}  \label{fhgsa2}.
\end{eqnarray}
\end{theorem}

Notice that when $g(\lambda)=\lambda^2/2$ and $f(\lambda)\equiv 0$, condition (\ref{condimenme1}) is called canonical assumption considered by de la Pe\~{n}a et al.\ \cite{D04,D07}. In particular, when $V_{i-1}$ is a constant and $g(\lambda)\equiv 0$, condition (\ref{condimenme1}) reduces to the condition  considered by Rio \cite{R13b}.

Next we show that Theorem \ref{th5} is very useful for obtaining the concentration inequalities for supermartingales.
Introducing the third moments of the supermartingale differences, we have the following Bernstein type inequalities.
\begin{corollary}\label{co9}Assume $\mathbf{E} (\xi _i^-)^3 < \infty$ for all $i \in [1, n]$.
 Denote by $
\langle\langle S \rangle\rangle_k=\sum_{i=1}^k \mathbf{E}((\xi_i^-)^3  |\mathcal{F}_{i-1})
$ for all $k \in [1, n].$
Then, for all $x, v, w>0$,
\begin{eqnarray}
&&\mathbf{P}\Big( S_k  \geq x,\ [S]_k \leq v^2\ \mbox{and}\  \langle\langle S \rangle\rangle_k\leq w  \ \mbox{for some}\ k\in[1,n]\Big)   \nonumber \\
&\leq&   \exp\left\{-\overline{\lambda} x + \frac12 \overline{\lambda}^2v^2 + \frac{1}3\overline{\lambda}^3w  \right\} \label{fhgst2} \\
&\leq&  B_1\left(x, \frac{w }{3v^2}, v\right)\label{dktgsdt1}\\
& \leq& B_2\left(x, \frac{w }{3v^2}, v\right),\label{dktgsdt2}
\end{eqnarray}
where $\overline{\lambda}=2x/(v^2+\sqrt{v^4+4w x})$.
\end{corollary}

Since $\langle\langle S \rangle\rangle_k  \leq \Upsilon(S_k)$, the inequalities (\ref{fhgst2}), (\ref{dktgsdt1}) and (\ref{dktgsdt2}) hold true when $\langle\langle S \rangle\rangle_k$ is replaced by $\Upsilon(S_k)$. To the best
of our knowledge, such inequalities have not been established for the sums of independent random variables.

Notice that (\ref{dktgsdt1}) and  (\ref{dktgsdt2}) are respectively the bounds of  Bennett and Bernstein.
Compared to the conditional Bernstein condition (\ref{Bernstein}), the condition of  Corollary \ref{co9} does not assume the existence of the moments of all orders.

For supermartingales with differences bounded from below, we still have the following Bernstein type inequality.
\begin{corollary}
\label{th7} Assume  $\xi_i \geq -1$ for all $i \in [1, n]$.
Then, for all $x, v>0$,
\begin{eqnarray}
 \mathbf{P}\Big( S_k  \geq x\ \mbox{and}\ [S]_k\leq v^2 \ \mbox{for some}\ k\in[1,n]\Big)
&\leq & \left(1+ \frac{x}{v^2}\right)^{v^2}e^{-x} \nonumber\\
&\leq& B_1\left(x, 1, v\right) \nonumber \\
&\leq& B_2\left(x, 1, v\right).\label{frie2}
\end{eqnarray}
\end{corollary}

Inequality (\ref{frie2}) is similar to Freedman's inequality (\ref{freedmanxf2}).
However, there are two differences between (\ref{frie2})  and  (\ref{freedmanxf2}). First, we assume $\xi_i$ bounded from below instead of $\xi_i$ bounded from above. Second,  the quadratic characteristic $\langle S\rangle_k$ in Freedman's inequality is replaced by the squared variation $[S]_k$ in our inequality (\ref{frie2}). Such inequality could be useful for
estimating the tail probabilities when the variances of $(\xi_i)$ do not exist.

Under the conditional Bernstein condition, we have
\begin{corollary}\label{co2}
 Assume, for a constant $\epsilon \in (0, \infty)$,
\begin{eqnarray}\label{br2}
 \mathbf{E}(\xi_{i}^{l}  | \mathcal{F}_{i-1})  \leq \frac12\, l!\, \epsilon^{l-2} \mathbf{E}(\xi_i^2|\mathcal{F}_{i-1}) \ \ a.s.\ \ \mbox{for all}\ \ l\geq 2\ \mbox{and all}\ i \in [1,n].
\end{eqnarray}
Then,
for all $  x, v> 0,$
\begin{eqnarray}
&& \mathbf{P}\Big( S_k \geq x\ \mbox{and}\ \langle S\rangle_{k}\leq v^2\ \mbox{for some}\ k \in [1,n] \Big) \nonumber\\
&\leq& B_{1,n}(x, \epsilon, v):=\exp \left\{- \overline{\lambda} x + n\log \left(1+ \frac{\overline{\lambda}^2 v^2  }{2n(1-\overline{\lambda}\epsilon)}\right)\right\}  \label{ie1g5}\\
&\leq& B_{1}(x, \epsilon, v) , \label{ie16}
\end{eqnarray}
where
\[
\overline{\lambda}= \frac{2x/v^2}{2x\epsilon/v^2+1+\sqrt{1+2x\epsilon/v^2}} \, \in (0, \epsilon^{-1}).
\]
\end{corollary}

Notice that $B_{1}(x, \epsilon, v) =\exp \left\{- \overline{\lambda} x +  \frac{\overline{\lambda}^2 v^2  }{2 (1-\overline{\lambda}\epsilon)} \right\}.$ In the independent case, inequality (\ref{ie16}) is known as Bennett's inequality \cite{B62}.
To highlight how the bound $B_{1,n}(x, \epsilon, v)$ improves Bennett's bound $B_{1}(x, \epsilon, v)$, we rewrite
\begin{eqnarray}
  B_{1,n}\left(x, \epsilon, v \right)  =  B_{1}(x, \epsilon, v) \exp\left\{-n \psi\left(\frac{\overline{\lambda}^2 v^2  }{2n(1-\overline{\lambda}\epsilon)}\right)\right\},\nonumber
\end{eqnarray}
where $\psi(t)=t-\log(1+t)$ is a nonnegative convex function in $t\geq 0$.   It is easy to see that, in the i.i.d.\ case with $v^2=n\sigma_1^2$ (or more generally when $\frac{\epsilon}{v}=\frac{\sigma_1}{\sqrt{n}}$ for a constant $\sigma_1>0$), we have
\begin{eqnarray}\label{fbnb}
 \ \ \ \ \ \ \ \ \  B_{1,n}\left(nx, \epsilon, \sqrt{n} \sigma_1 \right)  =  B_{1}\left(nx, \epsilon, \sqrt{n} \sigma_1 \right) \exp \left\{  -  n \, c_{x,\sigma_1,\epsilon} \right\},
\end{eqnarray}
where $c_{x,\sigma_1,\epsilon}=\psi\left(\frac{\overline{\lambda}^2 v^2  }{2 (1-\overline{\lambda}\epsilon)}\right)>0$ does not depend on $n$. Thus  Bennett's bound $B_{1}\left(nx, \epsilon, \sqrt{n} \sigma_1 \right)$ on tail probabilities $\mathbf{P}\left( S_n \geq nx \right)$ is strengthened by
adding a factor with exponential decay rate $\exp \left\{  -n  \, c_{x,\sigma_1,\epsilon} \right\}$ as $n\rightarrow \infty$.
Since the conditional Bernstein condition (\ref{Bernstein}) implies condition (\ref{br2}),  inequality  (\ref{ie1g5}) strengthen de la Pe\~{n}a's inequality (\ref{djknss}).

One calls $(\xi _i,\mathcal{F}_i)_{i=1,...,n}$ \emph{conditionally symmetric}, if $\mathbf{E}(\xi_i  >y |\mathcal{F}_{i-1})= \mathbf{E}(\xi_i < -y |\mathcal{F}_{i-1})$ for all $i \in [1,n]$ and for any $y \geq0$; see  Hitczenko \cite{H90},  de la Pe\~{n}a \cite{D99} and  Bercu and Touati \cite{BT08}. It is obvious that if $(\xi _i,\mathcal{F}_i)_{i=1,...,n}$ are conditionally symmetric, then, for any $y> 0$, $(\xi _i\mathbf1_{\{ |\xi_i|> y \}},\mathcal{F}_i)_{i=1,...,n}$  are also conditionally symmetric.
In particular, the conditionally symmetric martingale differences satisfy the canonical assumption $ \mathbf{E}(\exp\left\{\lambda \xi_i - \lambda^2\xi_i^2/2     \right\}|\mathcal{F}_{i-1})\leq 1$ for all $\lambda\geq 0$; see \cite{D99,D04,D07}. Thus, by Theorem \ref{th5} and optimizing on $\lambda$, inequality (\ref{fhgsa2}) implies de la Pe\~{n}a's inequality (\ref{dlpieq}).

The following result is a  Fuk-Nagaev  type inequality \cite{F73,N79} for martingales with conditionally symmetric differences. Its proof is based on a truncation argument on martingale differences.
\begin{corollary}\label{co1}
 Assume  that $(\xi _i,\mathcal{F}_i)_{i=1,...,n}$ are conditionally symmetric. Let
 \[
 V_k^2(y)=\sum_{i=1}^{k}\mathbf{E}(\xi_i^2\textbf{\emph{1}}_{\{|\xi_i|\leq y\}} |\mathcal{F}_{i-1}), \ \ \ k\in [1,n].
 \]
Then,  for all $x, y, v> 0$ and  $v^2\leq ny^2$,
\begin{eqnarray}
&& \mathbf{P}\Big(  S_k \geq x\ \mbox{and}\ V_k^2(y)\leq v^2\ \mbox{for some}\ k \in [1,n] \Big) \nonumber\\
&\leq&  \exp \left\{- \underline{\lambda} x + n\log \left(1+ \frac{ v^2  }{ n\,y^2}\left( \cosh(\underline{\lambda}\, y)  -1\right) \right)\right\}  +\, \mathbf{P}\left( \max_{1\leq i \leq n}  \xi_i  > y \right) \label{hdgh1} \\
&\leq&  \exp \left\{-  \overline{\lambda} x + \frac{ v^2  }{ y^2}\left( \cosh(\overline{\lambda}\, y)-1\right) \right\} +\, \mathbf{P}\left( \max_{1\leq i \leq n}  \xi_i > y \right), \label{hdgh2}
\end{eqnarray}
where
$$\underline{\lambda}=\frac{1}{y} \log\left( \frac{\frac{xy}{v^2}-\frac{x}{ny}+\sqrt{1+\frac{(xy)^2}{v^4} -2\frac{x^2}{nv^2} } }{1-\frac{x}{ny}} \,\right) \ \ \ \mbox{and}\ \ \ \ \overline{\lambda}=\frac{1}{y} \log \left( \frac{xy}{v^2} +\sqrt{1+\frac{(xy)^2}{v^4}  } \,\right).$$
\end{corollary}

Inequality (\ref{hdgh1}) is the best possible that can be obtained from the exponential Markov inequality $\mathbf{P}\left(  S_n \geq x \right) \leq  \inf_{\lambda\geq 0} \mathbf{E}e^{\lambda(S_n- x)}$
 under the present assumption. Indeed, if  $(\xi_i)_{i=1,...,n}$ are i.i.d. and satisfy the following  distribution
\begin{eqnarray}
\mathbf{P}(\xi_i=y  )= \mathbf{P}(\xi_i= -y )= \frac{v^2}{2ny^2}  \ \ \ \ \mbox{and} \ \ \ \mathbf{P}(\xi_i= 0 )= 1- \frac{v^2}{ ny^2},
\end{eqnarray}
then the bound (\ref{hdgh1}) equals to $\inf_{\lambda\geq 0} \mathbf{E}e^{\lambda(S_n- x)}$. In this sense, inequality (\ref{hdgh1}) is a version of Hoeffding's inequality (cf. (2.8) of \cite{Ho63}) for martingales with conditionally symmetric differences.

For martingales with bounded conditionally symmetric differences,  Sason \cite{S12} has obtained (\ref{hdgh1})
under the conditions $|\xi_i|\leq y$ and $\mathbf{E}(\xi_i^2 |\mathcal{F}_{i-1})\leq  v^2/n$. He has also obtained  (\ref{hdgh2}) under the assumption $|\xi_i|\leq y$. Thus (\ref{hdgh1}) improves and generalizes the Sason's inequalities under a more general condition.

For the martingales with square integrable differences,  several Nagaev type inequalities based on the truncation arguments on martingale differences can be found in Haeusler \cite{H84} and Courbot \cite{Co99}. For  optimal exponential convergence speed of such type bounds, we refer to  Lesigne and Voln\'{y} \cite{LV01} and Fan et al. \cite{F12,Fx1}.

Consider the case that the differences $ (\xi _i,\mathcal{F}_i)_{i=1,...,n}$ are sub-Gaussian. We have the following very general result.
\begin{theorem}
\label{th3} Assume that $V_{i-1}, i \in [1,n], $ are positive and $\mathcal{F}_{i-1}$-measureable random variables. Suppose   $ \mathbf{E}(e^{\lambda \xi_{i}}|\mathcal{F}_{i-1}) \leq  \exp \{ f(\lambda) V_{i-1}  \}$ for all $i \in [1,n]$ and for a positive function $ f(\lambda)$ for some $\lambda \in (0, \infty).$  Then, for all $ x, v > 0$,
\begin{eqnarray}
 \mathbf{P}\Big( S_{k} \geq x\ \mbox{and}\ \sum_{i=1}^{k} V_{i-1} \leq v^2 \ \mbox{for some}\ k \in [1,n]  \Big)
  \leq   \exp \Big\{ -\lambda x +f(\lambda)  v^2 \Big\} .\label{f1sd}
\end{eqnarray}
\end{theorem}

In the particular case where $v^2=\sum_{i=1}^{n} ||V_{i-1}||_{\infty}$ and  $f(\lambda)=\lambda^2/2$, Theorem \ref{th3} reduces to Theorem 4 of Fuk \cite{F73} after optimizing on $\lambda$. If $V_{i-1}=\mathbf{E}(\xi_i^2|\mathcal{F}_{i-1})$, Theorem \ref{th3} reduces to Theorem 4.2 of Khan \cite{K09}. Thus (\ref{f1sd}) can be regarded as
a generalization of the inequalities of  Fuk \cite{F73} and Khan \cite{K09}.

Using Theorem \ref{th3}, we extend  Azuma-Hoeffding's  inequality (cf.\ \cite{A67,Ho63}) to the case that the differences are only bounded from above.
\begin{corollary}
\label{co4} Assume that $U_{i-1}, i \in [1,n],$ are nonnegative and $\mathcal{F}_{i-1}$-measureable random variables.
Denote by
\begin{eqnarray}  \label{f1}
C_{i-1}^2 = \left\{ \begin{array}{ll}
\mathbf{E}(\xi_{i}^2|\mathcal{F}_{i-1}) , & \textrm{\ \ \ \ if  \ $ \mathbf{E}(\xi_{i}^2|\mathcal{F}_{i-1}) \geq U^2_{i-1}$ }, \\
\displaystyle\frac{1}{4}\left( U_{i-1}  + \frac{\mathbf{E}(\xi_{i}^2|\mathcal{F}_{i-1}) }{ U_{i-1} }\right)^2, & \textrm{\ \ \ \ otherwise}.
\end{array} \right.
\end{eqnarray}
If $\xi_{i} \leq U_{i-1}$  for all $i \in [1,n]$,
then, for all $\lambda >0$,
\begin{eqnarray}
 \mathbf{E}(e^{ \lambda \xi_i}|\mathcal{F}_{i-1}) \leq \exp\left\{ \frac{\lambda^2}{2} C_{i-1}^2 \right\} ,\label{gnlm}
\end{eqnarray}
and,  for all $ x, v >0$,
\begin{eqnarray}\label{f3}
 \mathbf{P}\bigg(  S_k \geq x\  \mbox{and}\  \sum_{i=1}^k C_{i-1}^2 \leq v^2\ \mbox{for some}\ k \in [1, n]  \bigg)  \leq   \exp\left\{- \frac{ x^2}{  2\, v^2 }   \right\}.
\end{eqnarray}
In particular, if $\mathbf{E}(\xi_{i}^2|\mathcal{F}_{i-1}) \geq U^2_{i-1}$ for all $i\in [1,n]$, then, for all $x,v \geq0$,
\begin{eqnarray}\label{f3saf}
 \mathbf{P}\left( \max_{1\leq k \leq n} S_k \geq x\  \mbox{and}\  \langle S\rangle_n \leq v^2\right)  \leq   \exp\left\{- \frac{ x^2}{  2\, v^2 }   \right\}.
\end{eqnarray}
\end{corollary}

Notice that if $(\xi_i)_{i=1,...,n}$ are independent and satisfy
the conditions $\xi_i \leq c_i$ and $\mathbf{E}\xi_i^2\geq c_i^2$ for some constants $(c_i)_{i=1,...,n}$, then (\ref{f3saf}) is a gaussian bound with $v^2= \sum_{i=1}^n\mathbf{E}\xi_i^2$.
It is obvious that the Rademacher random variables satisfy this assumption.

For martingale differences $(\xi _i,\mathcal{F}_i)_{i=1,...,n}$, inequality
(\ref{f3}) generalizes the following inequality due to van de Geer (cf.\
Theorem 2.5 of \cite{V02}): if $L_{i-1} \leq \xi_{i} \leq U_{i-1}$ for some
$\mathcal{F}_{i-1}$-measureable random variables $L_{i-1}$ and $U_{i-1}$,
then, for all $x,v > 0$,
\begin{eqnarray}\label{fgr}
 \mathbf{P}\bigg(  S_k \geq x\  \mbox{and}\ \frac14 \sum_{i=1}^k \left( U_{i-1} -L_{i-1}\right)^2 \leq v^2\ \mbox{for some}\ k \in [1, n]  \bigg)  \leq   \exp\left\{- \frac{ x^2}{ 2\, v^2 }   \right\}.
\end{eqnarray}
Indeed, since
\begin{eqnarray*}
\mathbf{E}(\xi_{i}^2|\mathcal{F}_{i-1})= \mathbf{E}((\xi_{i}-L_{i-1}) \xi_{i}|\mathcal{F}_{i-1})  \leq  \mathbf{E}((\xi_{i}-L_{i-1}) U_{i-1}|\mathcal{F}_{i-1}) \leq -L_{i-1}U_{i-1}\ ,
\end{eqnarray*}
we have
\begin{eqnarray}
\sum_{i=1}^k C_{i-1}^2  \leq   \sum_{i=1}^k \frac{1}{4}\left( U_{i-1}  + \frac{\mathbf{E}(\xi_{i}^2|\mathcal{F}_{i-1}) }{ U_{i-1} }\right)^2   \leq   \frac{1}{4} \sum_{i=1}^k \left( U_{i-1} -L_{i-1}\right)^2 \label{fsav}
\end{eqnarray}
and
\[
\left\{ \frac{1}{4} \sum_{i=1}^k \left( U_{i-1} -L_{i-1}\right)^2 \leq v^2 \right\} \subseteq \left\{  \sum_{i=1}^k C_{i-1}^2 \leq v^2 \right\},
\]
which together with (\ref{f3}) implies (\ref{fgr}).

Under the assumption of Corollary \ref{co4}, Pinelis \cite{P06,P06b} (see also Bentkus \cite{Be03}) proved the following inequality, for all $x>0$,
\begin{eqnarray}
 \mathbf{P}\bigg( \max_{1\leq k \leq n} S_k \geq x
 \bigg)  &\leq& c\left(1-\Phi(\frac{x}{\hat{v}})\right) \label{pio}\\
  &=& O\left( \frac{1}{1+ x/  \hat{v}  }  \exp\left\{- \frac{   x^2}{ 2\,  \hat{v}^2 }   \right\} \right),\ \ \ \ x\rightarrow \infty, \nonumber
\end{eqnarray}
where $c$ is an absolute constant and $$\hat{v}^2= \sum_{i=1}^n \left|\left|\frac{1}{4}\left( U_{i-1}  + \frac{\mathbf{E}(\xi_{i}^2|\mathcal{F}_{i-1}) }{ U_{i-1} }\right)^2\right|\right|_{\infty}.$$
Notice that $\hat{v}^2 \geq \left|\left|\sum_{i=1}^n C_{i-1}^2\right|\right|_{\infty}.$
If $\hat{v}^2 = \left|\left|\sum_{i=1}^n C_{i-1}^2\right|\right|_{\infty}$, Pinelis' inequality (\ref{pio})
is better than ours (\ref{f3}) by adding a factor $\frac{O(1)}{1+ x/  \hat{v} }$. Otherwise $\hat{v}^2 > \left|\left|\sum_{i=1}^n C_{i-1}^2\right|\right|_{\infty}$,   our inequality (\ref{f3}) improves  Pinelis' inequality (\ref{pio}) by adding an exponential decay factor of order
\begin{eqnarray}\label{mspi}
 \left(1+ \frac{x}{\hat{v}} \right) \exp\left\{ -\frac{x^2}{2}\delta \right\}, \ \ \ \ \ x\rightarrow \infty,
\end{eqnarray}
where $$\delta=\frac{ \hat{v}^2-\left|\left|\sum_{i=1}^n C_{i-1}^2\right|\right|_{\infty}}{ \hat{v}^2 \left|\left|\sum_{i=1}^n C_{i-1}^2\right|\right|_{\infty} } > 0.$$
To illustrate this factor, consider the following example. For a  much more significant improvement, we refer to  Remark \ref{fsfdf}.

\emph{Example 1}: Assume that $(\varepsilon_{i})_{i=1,...,n}$ is a sequence of Rademacher random variables, and that $\mathcal{N}$ is a random variable independent of $(\varepsilon_{i})_{i=1,...,n}$. Set
\begin{eqnarray*}
\xi_i = \left(\frac{ \varepsilon_{i}}{\sqrt{n}}   \sin \mathcal{N}\,\right) \mathbf{1}_{\{i\, \textrm{ is odd}\}} +  \left( \frac{ \varepsilon_{i}}{\sqrt{n}}  \cos  \mathcal{N}\, \right) \mathbf{1}_{\{i\, \textrm{ is even}\}},
\end{eqnarray*}
 $\mathcal{F}_0=\sigma\{\mathcal{N} \}$ and $\mathcal{F}_i=\sigma\{\mathcal{N} ,$ $ \varepsilon_{j}, 1\leq j \leq i \}$. So we have
\begin{eqnarray*}
\xi_i \leq U_{i} :=  \left(\frac{ 1}{\sqrt{n}}  | \sin \mathcal{N} |\,\right) \mathbf{1}_{\{i\, \textrm{ is odd}\}}   +   \left(\frac{ 1}{\sqrt{n}}  | \cos \mathcal{N} |\,\right)\, \mathbf{1}_{\{i\, \textrm{ is even}\}}
\end{eqnarray*}
and
\begin{eqnarray*}
\sum_{i=1}^n C_{i-1}^2 = \langle S\rangle_n=\sum_{i=1}^n \left( \frac{ \sin^2 \mathcal{N}}{ n } \, \mathbf{1}_{\{i\, \textrm{ is odd}\}} +   \frac{ \cos^2 \mathcal{N}}{ n }  \, \mathbf{1}_{\{i\, \textrm{ is even}\}} \right).
\end{eqnarray*}
Hence, for any even number $n$, it is easy to see that  $\hat{v}^2=1 > \frac{1}{2}= \sum_{i=1}^n C_{i-1}^2=\langle S\rangle_n $.
Then Pinelis'  inequality (\ref{pio}) shows that:
$$\mathbf{P}\left( \max_{1\leq k \leq n} S_k \geq x
 \right)  = O\bigg( \frac{1}{1+ x   }  \exp\left\{- \frac{ \,  x^2}{ 2 }   \right\} \bigg),\ \ \ \ x \rightarrow \infty, $$
while our inequality (\ref{f3}) implies that:
$$\mathbf{P}\Big( \max_{1\leq k \leq n} S_k \geq x
 \Big)  \leq   \exp\Big\{- x^2  \Big\} .$$
Thus our inequality (\ref{f3}) improves Pinelis' inequality (\ref{pio}) by adding a factor with the exponential decay rate $(1+ x)\exp\left\{- \frac{x^2}{ 2} \right\} .$

\begin{remark} \label{eton} Corollary \ref{co4} implies a simple proof of the following self-normalized deviation inequality.  Assume  that $(\xi _i)_{i=1,...,n}$ are independent and symmetric.
Then, for all $x > 0$,
\begin{eqnarray}
 \mathbf{P}\left( \max_{1\leq k \leq n} \frac{S_k}{\sqrt{[S]_n}}  \geq x \right)
 \leq   \exp\left\{- \frac{x^2}{ 2  }   \right\}, \label{fdgfgfdgh}
\end{eqnarray}
where by convention $\frac00=0.$ A similar result can be found in  Hitczenko \cite{H90}. Hitczenko  has obtained the same upper bound on tail probabilities $\mathbf{P}\left( S_n   \geq x ||\sqrt{[S]_n} ||_{\infty}\right).$
For more precise results, we refer to  Wang and Jing \cite{WJ99}. In particular, the Cram\'{e}r type large deviations have been established  by  Jing, Shao and Wang \cite{JSW03} without assuming that $(\xi _i)_{i=1,...,n}$ are symmetric (or $(\xi _i)_{i=1,...,n}$ have exponential moments).
\end{remark}

\section{Applications to statical estimation}\label{sec2.5}
The exponential concentration inequalities for martingales certainly have many
applications. \mbox{McDiarmid} \cite{M} and Rio \cite{R13a} applied such type
inequalities to estimate the concentration of separately Lipschhitz functions.
Van de Geer \cite{V95} found that such inequalities can be used for maximum
likelihood estimation for counting processes. Liu and Watbled \cite{Liu09a}
considered the free energy of directed polymers in a random environment via martingale inequalities.
Dedecker and Fan \cite{DF12} gave an application of these inequalities to the
Wasserstein distance between the empirical measure and the invariant
distribution. We refer to Bercu \cite{B08} for more interesting applications
of the concentration inequalities for martingales.

In the sequel, we discuss how to apply our results to linear regression models, autoregressive processes and branching processes. We find these models in Liptser and Spokoiny \cite{Ls01}  and Bercu and Touati \cite{BT08}.

\textbf{1. Linear regression models.} Consider the stochastic linear regression models given, for all $k \in [1, n],$ by
\begin{equation}\label{ine29}
X_{k}=\theta \phi_k + \varepsilon_{k}
\end{equation}
where $X_k, \phi_k$ and $\varepsilon_{k}$ are the observations, the regression variables and the driven noises, respectively. We assume that $(\phi_k)$ is a sequence of independent  random variables. We also assume that $(\varepsilon_k)$ is a sequence of independent  and identically distributed (i.i.d.) random variables, with mean zero and variation $\sigma^2>0.$ Moreover, we suppose that $(\phi_k)$ and  $(\varepsilon_k)$ are independent. Our interest is to estimate the unknown parameter $\theta.$ The well-known  least-squares estimator $\theta_n$ is given below
\begin{equation}\label{ine30}
\theta_n = \frac{\sum_{k=1}^n \phi_{k} X_k}{\sum_{k=1}^n \phi_{k}^2}.
\end{equation}
When $(\phi_k)$ and  $(\varepsilon_k)$ are sub-Gaussian,  exponential inequalities on the convergence of $\theta_n -\theta$ have been established by Bercu and Touati \cite{BT08}.   When $(\varepsilon_k)$ are the normal random variables,    Liptser and Spokoiny  \cite{Ls01} have established the following estimation: for all $x\geq 1,$
\begin{eqnarray}\label{ls21}
 \mathbf{P}\left( \pm \, (\theta_n -\theta)\sqrt{ \Sigma _{k=1}^n \phi_{k}^2} \geq x \sigma  \right) \  \leq\    \sqrt{\frac2 \pi} \, \frac{1}{x}   \, \exp\bigg\{ - \frac{x^2}{2  }  \bigg\}.
\end{eqnarray}
Here, we would like to give a generalization of this inequality. Consider the case that the random variables $(\varepsilon_k)$ satisfy the Bernstein condition.
\begin{theorem} \label{thlin}
Assume  $ |\phi_{k}|/\sqrt{\sum_{k=1}^n \phi_{k}^2} \leq \epsilon_1$ and
 \[
|\mathbf{E} \varepsilon_{i}^{k}   | \leq \frac12 k!\epsilon_2^{k-2} \mathbf{E} \varepsilon_{i}^2  ,\ \ \   \ \ \textrm{for all}\ k\geq 2\ \textrm{and all}\  i \in [1, n],
\]
for  two positive numbers $\epsilon_1$ and $\epsilon_2$. Let $\epsilon = \epsilon_1 \epsilon_2/\sigma.$
Then, for all $x \geq 0$,
\begin{eqnarray}\label{th5ineq}
 \mathbf{P}\left( \pm \, (\theta_n -\theta)\sqrt{ \Sigma _{k=1}^n \phi_{k}^2} \geq x \sigma  \right) \  \leq\   B_{1,n}(x, \epsilon, 1) \ \leq \ \exp\bigg\{ - \frac{x^2}{2(1+ x\epsilon) }  \bigg\}.
\end{eqnarray}
\end{theorem}

Since $ |\phi_{k}|/\sqrt{\sum_{k=1}^n \phi_{k}^2} \leq 1,$ the condition imposed on $(\phi_{k})$ of Theorem \ref{thlin} can be dropped  by taking $\epsilon_1 =1$.
It is interesting to see that by taking  $\epsilon_1 =1,$  bound (\ref{th5ineq}) does not depend on the distribution of the regression variables $(\phi_{k})$. This is a big advantage in practice.

If $a\leq |\phi_{k}| \leq b $ for two positive constants $a$ and $b$, then the condition of Theorem \ref{th5} is  satisfied with $\epsilon_1 = \frac{b\, }{a \sqrt{n  } }.$  Indeed, it is easy to see that
\[
\frac{|\phi_{k}|}{\sqrt{\sum_{i=1}^n \phi_{i}^2} }  \leq   \frac{b}{\sqrt{n a^2} }  =  \epsilon_1.
\]
In this case, bound (\ref{th5ineq})  behaviors like $\exp\{-x^2/2\}$ when $x=o( \sqrt{n})$ as $n\rightarrow \infty$. When $x$ is large, bound (\ref{th5ineq}) behaviors like $\exp\{-x  \}$.

If $(\varepsilon_k)$ are  bounded from above, we have the following sub-Gaussian  tail bound  from Corollary \ref{co4}.
\begin{theorem}\label{dssaf}
If $ \varepsilon_k \leq \epsilon$ for all $k \in [1, n],$ then, for all $x \geq 0$,
\begin{eqnarray} \label{tsfsdfd}
 \mathbf{P}\left(   (\theta_n -\theta)\sqrt{ \Sigma _{k=1}^n \phi_{k}^2} \geq x \sigma  \right)  \ \leq \      \exp \bigg\{ -\frac{x^2 }{2 C_n } \bigg\},
\end{eqnarray}
where
\[
C_n = \frac14  \Big( \frac{ \epsilon}{ \sigma}  + \frac{\sigma }{\epsilon} \Big)^2  .
\]
In particular, if $ |\varepsilon_k| \leq \epsilon,$  bound (\ref{tsfsdfd})  holds true on the tail probabilities
 $$\mathbf{P}\left( \pm \, (\theta_n -\theta)\sqrt{ \Sigma _{k=1}^n \phi_{k}^2} \geq x \sigma  \right). $$
\end{theorem}

\begin{remark}\label{fsfdf} If $ |\varepsilon_k| \leq \epsilon,$  we can obtain some similar bounds by using van de Geer's inequality (\ref{fgr}) or Pinelis' inequality (\ref{pio}).  However, those bounds are less tight than  (\ref{tsfsdfd}). Indeed, by van de Geer's inequality,
  we can obtain the bound (\ref{tsfsdfd}) with a larger
  $C_n=(\epsilon/\sigma)^2.$ If we make use of Pinelis' inequality (or
  Bentkus' inequality \cite{Be03}), the bound will be as large as $$
  \frac{O(1)}{x} \exp \bigg\{ -\frac{x^2 }{2 n C_n } \bigg\}. $$
\end{remark}

Next, consider the tail probabilities of $(\theta_n -\theta) \sum_{k=1}^n \phi_{k}^2.$ It seems that our inequalities fit well
to such type estimations.
\begin{theorem}\label{thds}
Assume that there exist  $ \alpha \in (1, 2]$  and $c> 0$ such that
 \[
 \mathbf{E}e^{ \lambda \varepsilon_{i} }   \leq e^{c  |\lambda|^\alpha}  \ \ \  \ \ \textrm{for all }   i\in [1, n] \textrm{ and all } \lambda \in \mathbf{R} .
\]
Then, for all $x,   v \geq 0$,
\begin{eqnarray}\label{sfcscs}
 \mathbf{P}\Big(   \pm (\theta_n -\theta)  \sum_{k=1}^n \phi_{k}^2  \geq x  \  \textrm{and} \   \sum_{k=1}^n   |\phi_{k}|^\alpha \leq v^\alpha  \Big)  &\leq&    \exp\left\{-C(\alpha) \left(\frac{x}{v} \right)^\frac{\alpha}{\alpha-1} \right\},
\end{eqnarray}
where $$C(\alpha)=(c\, \alpha )^{\frac1{1-\alpha}} \left(1- \alpha^{-1} \right).$$
\end{theorem}

When the condition of Theorem \ref{thds} is verified with $\alpha=2,$ then $(\varepsilon_i)$ are known as sub-Gaussian random variables. It is known that the bounded random variables and the normal random variables are all sub-Gaussian random variables.
In particular, if $(\varepsilon_i)$ are the standard normal random variables, then bound (\ref{sfcscs}) is valid with $\alpha=2$ and $c=C(2)=1/2$.

\textbf{2. Autoregressive processes.} The model of autoregressive can be stated as follows: for all $k \in [1,  n],$
\begin{equation}\label{scfso1}
X_{k} = \theta X_{k-1} +  \varepsilon_k \,  ,
\end{equation}
where $(X_k)$ and $(\varepsilon_k)$ are the observations and driven noises, respectively. We assume  that $(\varepsilon_k)$ is a sequence of i.i.d.\ centered random variables with  variation $\sigma^2>0.$ The process is said to be stable if $|\theta|\leq 1,$ unstable if $|\theta|=1$ and explosive if $|\theta|>1.$ We can estimate the unknown parameter $\theta$ by the least-squares  estimator given by, for all $n\geq 1$,
\begin{equation}\label{scfso2}
\theta'_n  = \frac{\sum_{k=1}^n X_{k}X_{k-1}}{\sum_{k=1}^n X_{k-1}^2}  .
\end{equation}
When $X_0$ and $(\varepsilon_k)$ are the normal random variables, the convergence rate of $\theta_n '-\theta$  has been established by Bercu and Touati \cite{BT08}.  Here, we would like to give an almost sure convergence rate of $(\theta_n' -\theta) \sum_{k=1}^n X_{k-1}^2 .$

By an argument similar to that of Theorem \ref{thds}, we have the following result.
\begin{theorem}
Assume the condition of Theorem \ref{thds}.
Then bound (\ref{sfcscs}) holds true on the tail probabilities
\begin{eqnarray}\label{sfcsc2s}
 \mathbf{P}\Big(   \pm (\theta_n' -\theta)  \sum_{k=1}^n  X_{k-1}^2  \geq x  \  \textrm{and} \   \sum_{k=1}^n   | X_{k-1}|^\alpha \leq v^\alpha  \Big).
\end{eqnarray}
\end{theorem}

If $(\varepsilon_{i})$ are bounded, then we have
\begin{theorem}\label{th6}
Assume $|\varepsilon_{i}| \leq \epsilon$ for all $i \in [1,   n].$
Then, for all $x, v > 0$,
\begin{eqnarray}\label{tsfsdfdog}
 \mathbf{P}\left(\pm \,  (\theta_n '-\theta)  \sum _{k=1}^n X_{k-1}^2  \geq x \  \textrm{and} \ L_n  \leq v^2 \right) \  \leq\    \exp \bigg\{ -\frac{x^2 }{2\, v^2} \bigg\},
\end{eqnarray}
where
\[
L_n = \frac14  \Big( \epsilon + \frac{\sigma^2}{\epsilon} \Big)^2 \sum _{k=1}^n X_{k-1}^2.
\]
\end{theorem}
\begin{remark}\label{fsfdf1} We can obtain some similar bounds by using Corollary 2.3 or van de Geer's inequality.  However, those bounds are less tight than  (\ref{tsfsdfdog}). For instance, by van de Geer's inequality,
  we can obtain the bound (\ref{tsfsdfdog}) with a larger $L_n= \epsilon^2
  \sum _{k=1}^n X_{k-1}^2.$
\end{remark}

\textbf{3. Branching processes.}  Consider the Galton-Watson process stating from $X_0=1$ and given, for all $n\geq 1$, by
\[
X_n = \sum_{k=1}^{X_{n-1}} Y_{n, k} \, ,
\]
where $(Y_{n, k})$ is a sequence i.i.d.\ and nonnegative integer-valued random variables. The distribution of $(Y_{n, k}),$ with finite mean $m$ and variance $\sigma^2,$ is commonly called the offspring or reproduction distribution. We are interested in the estimation of the offspring  mean $m.$ The Lotka-Nagaev estimator is given by
\[
m_n=  \frac{X_{n}}{X_{n-1}} \,.
\]
Assume $X_{n}> 0$ a.s.\ such that the Lotka-Nagaev estimator $m_n$ is always well defined. Our goal is to establish exponential inequalities
for  $m_n.$ Denote by
$$ \xi_{n, k} = Y_{n, k} -m .$$
Then
\[
(m_n  -m ) X_{n-1} = X_n -mX_{n-1} = \sum_{k=1}^{X_{n-1}}\xi_{n, k} \, .
\]
Thus $(m_n  -m ) X_{n-1}$ is a sum of independent random variables by given $X_{n-1}$.
By Corollary \ref{co2}, we easily obtain the following exponential inequalities.
\begin{theorem}
Assume, for a constant $\epsilon \in (0, \infty)$,
\begin{eqnarray*}
 |\mathbf{E} \xi_{n, k}^{l}  |  \leq \frac12\, l!\, \epsilon^{l-2} \mathbf{E} \xi_{n, k}^2  \ \ \ \ \mbox{for all}\ \ l\geq 2\ \mbox{and all}\ k \in [1,X_{n-1}].
\end{eqnarray*}
Then, for all $x, v > 0$, it holds
\begin{eqnarray*}
 \mathbf{P}\Big(  |m_n  -m | X_{n-1} \geq x  \  \textrm{and} \ X_{n-1}\sigma^2   \leq v^2 \, \Big| \, X_{n-1} \Big) & \leq&    2 \, B_{1,n}(x, \epsilon,   v) \\
  &\leq& 2 \exp\bigg\{ - \frac{x^2}{2(  v^2+ x\epsilon) }  \bigg\}.
\end{eqnarray*}
In particular, it implies that, for all $x  > 0,$
\begin{eqnarray*}
 \mathbf{P}\Big(  |m_n  -m |  \geq  x  \Big) \    \leq \  2\, \mathbf{E}\bigg( \exp \bigg\{ - X_{n-1} \frac{x^2 }{2\, (\sigma^2+ x\epsilon) } \bigg\} \, \bigg).
\end{eqnarray*}
\end{theorem}

Since $\xi_{n, k} \geq -m$,  we have the following one side sub-Gaussian bound by Corollary
\ref{co4}. This bound cannot be obtained from  Azuma-Hoefding's inequality.
\begin{theorem} For all $x, v > 0$,  it holds
\begin{eqnarray*}
 \mathbf{P}\Big(  (m_n  -m ) X_{n-1} \leq - x, M  X_{n-1} \leq v^2 \, \Big| \, X_{n-1} \Big) \    \leq \   \exp \bigg\{ -\frac{x^2 }{2\, v^2 } \bigg\} \, ,
\end{eqnarray*}
where
\begin{eqnarray*}  \label{f1}
M   = \left\{ \begin{array}{ll}
\sigma^2 , & \textrm{\ \ \ \ \ if $\sigma \geq m$, } \\
\displaystyle\frac{1}{4}\left( m  + \frac{\sigma^2}{ m }\right)^2, & \textrm{\ \ \ \ \ if $\sigma < m$.}
\end{array} \right.
\end{eqnarray*}
In particular, it implies that, for all $x  > 0,$
\begin{eqnarray*}
 \mathbf{P}\Big(  (m_n  -m )   \leq - x  \Big) \    \leq \   \mathbf{E}\bigg( \exp \bigg\{ - X_{n-1} \frac{x^2 }{2\, M } \bigg\} \, \bigg).
\end{eqnarray*}
\end{theorem}

More generale estimations on the tail probabilities $\mathbf{P}\left(  |m_n  -m | \geq  x  \right),$ we refer to Bercu and Touati \cite{BT08}. In particular,  Bercu and Touati have established the  Bernstein bounds associated with the cumulant generating function of $\xi_{n, k}.$

\section{ Proof of Theorem \ref{th5} }\label{sec6}
Suppose $ \mathbf{E}(\exp\left\{\lambda \xi_i - g(\lambda)\xi_i^2     \right\}|\mathcal{F}_{i-1}) \leq  1+ f(\lambda)V_{i-1}$ for a constant $ \lambda \in (0, \infty)$ and all $i \in [1, n].$ Define the \emph{exponential multiplicative
martingale} $Z(\lambda )=(Z_k(\lambda ),\mathcal{F}_k)_{k=0,...,n},$ where
\[
   Z_k(\lambda )=\prod_{i=1}^{  k}\frac{\exp\left\{\lambda \xi_i -g (\lambda)\xi_i^2     \right\}}{\mathbf{E}\left(\exp\left\{\lambda \xi_i -g(\lambda)\xi_i^2     \right\} |
\mathcal{F}_{i-1} \right)},  \quad \quad  \quad Z_0(\lambda )=1. \label{C-1}
\]
If $T$ is a stopping time, then $Z_{T\wedge k}(\lambda )$ is also a martingale, where
\[
Z_{T\wedge k}(\lambda )=\prod_{i=1}^{T\wedge k}\frac{\exp\left\{\lambda \xi_i -g (\lambda)\xi_i^2     \right\}}{\mathbf{E}\left(\exp\left\{\lambda \xi_i -g(\lambda)\xi_i^2     \right\} |
\mathcal{F}_{i-1} \right)}, \quad \quad Z_0(\lambda )=1.
\]
Thus,  the random variable $Z_{T\wedge k}(\lambda ) $ is a
probability density on $(\Omega ,\mathcal{F},\mathbf{P})$, i.e.
$$ \int Z_{T\wedge k}(\lambda)  d \mathbf{P} = \mathbf{E}(Z_{T\wedge k}(\lambda))=1.$$
 Define the \emph{conjugate probability measure}
\begin{equation}
d\mathbf{P}_\lambda =Z_{T\wedge n}(\lambda )d\mathbf{P}.  \label{chmeasure3}
\end{equation}
Denote $\mathbf{E}_{\lambda}$ the expectation with
respect to $\mathbf{P}_{\lambda}.$

\vspace{0.3cm}

\noindent\emph{Proof of Theorem \ref{th5}.} For any $x,  v, w>0$, define the stopping time
\[
T(x,v,w)=\min\left\{k\in [1, n]: S_{k} \geq x,\ [S]_{k} \leq v^2 \ \mbox{and}\ \sum_{i=1}^{k} V_{i-1} \leq w  \right\},
\]
with the convention that $\min{\emptyset}=0$. Then
 \[
\textbf{1}_{\{ S_{k} \geq x,\ [S]_{k} \leq v^2 \ \mbox{and}\ \sum_{i=1}^{k} V_{i-1} \leq w\ \mbox{for some}\ k\in[1,n] \}} = \sum_{k=1}^{n}  \textbf{1}_{\{ T(x,v,w)=k\}}.
\]
By the change of measure (\ref{chmeasure3}), we deduce that, for all $x,\lambda, v,w>0$,
\begin{eqnarray}
  && \mathbf{P}\left( S_{k} \geq x,\ [S]_{k} \leq v^2 \ \mbox{and}\ \sum_{i=1}^{k} V_{i-1} \leq w\ \mbox{for some}\ k\in[1,n] \right) \nonumber\\
 &=& \mathbf{E}_{\lambda}\Big( Z_{T\wedge n}(\lambda)^{-1}\textbf{1}_{\{S_{k} \geq x,\ [S]_{k} \leq v^2 \ \mbox{and}\ \sum_{i=1}^{k} V_{i-1} \leq w\ \mbox{for some}\ k\in[1,n]\}}\Big) \nonumber \\
 &=& \sum_{k=1}^{n}\mathbf{E}_{\lambda}\Big(  \exp\{-\lambda S_{k}+g(\lambda) [S]_{k} + \Xi_{k}(\lambda) \} \textbf{1}_{\{T(x,v,w)=k\}}\Big)  , \label{ghnda}
\end{eqnarray}
where
\[
\Xi_{k}(\lambda)= \sum_{i=1}^k \log \mathbf{E} \left(\exp\left\{\lambda \xi_i -g(\lambda)\xi_i^2 \right\} |\mathcal{F}_{i-1}\right).
\]
Using Jensen's inequality and the condition $ \mathbf{E}(\exp\left\{\lambda \xi_i - g(\lambda)\xi_i^2     \right\}|\mathcal{F}_{i-1}) \leq  1+ f(\lambda)V_{i-1}$, we have
\begin{eqnarray}
\Xi_{k}(\lambda) &\leq& k \log\left(   \frac{1}{k} \sum_{i=1}^k \mathbf{E} \left(\exp\left\{\lambda \xi_i -g(\lambda)\xi_i^2 \right\} |\mathcal{F}_{i-1}\right) \right) \nonumber \\
&\leq& k \log\left( 1+ \frac{1}{k} f(\lambda)\sum_{i=1}^k V_{i-1} \right).
\end{eqnarray}
Thus (\ref{ghnda}) implies that,
for all $x,\lambda, v, w>0$,
\begin{eqnarray}
  && \mathbf{P}\left( S_{k} \geq x,\ [S]_{k} \leq v^2 \ \mbox{and}\ \sum_{i=1}^{k} V_{i-1} \leq w\ \mbox{for some}\ k\in[1,n] \right) \nonumber\\
 &\leq& \sum_{k=1}^{n}\mathbf{E}_{\lambda} \Bigg( \exp \left\{-\lambda S_k +g(\lambda)[S]_{k}  + k \log\left( 1+ \frac{1}{k} f(\lambda)\sum_{i=1}^k V_{i-1} \right) \right\} \textbf{1}_{\{T(x,v,w)=k\}}\Bigg) .\ \ \ \
\end{eqnarray}
By the fact $S_{k}\geq x$, $[S]_{k}\leq v^2$ and $\sum_{i=1}^k V_{i-1}\leq  w$  on the set  $\{T(x,v,w)=k\}$, we find that,
for all $x,\lambda, v, w>0$,
\begin{eqnarray}
  && \mathbf{P}\left( S_{k} \geq x,\ [S]_{k} \leq v^2 \ \mbox{and}\ \sum_{i=1}^{k} V_{i-1} \leq w\ \mbox{for some}\ k\in[1,n] \right) \nonumber\\
 &\leq&   \exp \left\{-\lambda x +g(\lambda)v^2  + k \log\left( 1+ \frac{1}{k} f(\lambda)w \right) \right\}  \mathbf{E}_{\lambda} \Big(\sum_{k=1}^{n}\textbf{1}_{\{T(x,v,w)=k\}}\Big)\nonumber\\
 &\leq&   \exp \left\{-\lambda x +g(\lambda)v^2  + n \log\left( 1+ \frac{1}{n} f(\lambda)w \right) \right\}  \label{sthv1}\\
&\leq&   \exp \left\{-\lambda x +g(\lambda)v^2  + f(\lambda)w \right\}  \label{sthv2}.
\end{eqnarray}
This gives the desired inequalities (\ref{fhgsa1}) and (\ref{fhgsa2}), and completes the proof of Theorem \ref{th5}.

\vspace{0.3cm}

\noindent\emph{Proof of Corollary \ref{co9}.} To prove Corollary \ref{co9}, we should use the following basic inequality:
\[
\exp\left\{x-\frac{1}{2}x^2\right\} \leq 1 + x +\frac{1}{3}(x^-)^3 ,\ \ \ \ x \in \mathbf{R}.
\]
By the last inequality, it follows that, for all $\lambda>0$,
\begin{eqnarray}
 \mathbf{E}\left(\exp\left\{\lambda \xi_i -\frac{1}{2} (\lambda\xi_i)^2     \right\} \bigg|\mathcal{F}_{i-1} \right) \leq  1+ \frac{1}{3}\lambda^3 \mathbf{E}\left( (\xi_i^-  )^3 |\mathcal{F}_{i-1} \right).
\end{eqnarray}
Applying the inequalities (\ref{fhgsa1}) and (\ref{fhgsa2}) with
$g(\lambda)=\frac{\lambda^2}{2}$, $f(\lambda)=\frac{\lambda^3}{3}$ and $V_{i-1}=\mathbf{E}\left( (\xi_i^-  )^3 |\mathcal{F}_{i-1} \right),$
 we get  (\ref{fhgst2})  by noting the fact that
\begin{eqnarray}\label{ineq385}
 \inf_{\lambda>0}\exp\left\{-\lambda x + \frac12  \lambda^2v^2 + \frac{1}3\lambda^3w \right\} =  \exp\left\{-\overline{\lambda} x + \frac12  \overline{\lambda}^2v^2 + \frac{1}3\overline{\lambda}^3w \right\},
\end{eqnarray}
where $\overline{\lambda}=2x/(v^2+\sqrt{v^4+4wx}).$ By a simple calculation, we find that,  for all $v, w >0$ and all $0<\lambda < \frac{3v^2}{w}$,
\begin{eqnarray}
\frac12  \lambda^2v^2 + \frac{1}3\lambda^3w  \leq  \frac{\lambda^2v^2}{2( 1-\frac{\lambda w}{3 v^2})}.
\end{eqnarray}
Thus, for all $x,  v, w >0$,
\begin{eqnarray}
  \exp\left\{-\overline{\lambda} x + \frac12  \overline{\lambda}^2v^2 + \frac{1}3\overline{\lambda}^3w \right\} &\leq&   \inf_{0< \lambda < \frac{3v^2}{w}}\exp\left\{-\lambda x + \frac{\lambda^2v^2}{2( 1-\frac{\lambda w}{3 v^2})} \right\} \nonumber\\
 &=& B_1\left(x, \frac{w}{3v^2}, v\right)\nonumber\\
 &\leq& B_2\left(x, \frac{w}{3v^2}, v\right).\nonumber
\end{eqnarray}
Combining this inequality with (\ref{fhgst2}), we obtain the desired inequalities (\ref{dktgsdt1}) and (\ref{dktgsdt2})  of  the corollary.

To prove Corollary \ref{th7}, we need the following lemma. 
\begin{lemma}
\label{lemma17} If $\xi$ is a random variable such that $\xi \geq -1$ and $\mathbf{E}\xi
\leq 0$, then, for all $\lambda \in  [0, 1),$
\begin{eqnarray*}
 \mathbf{E}\Big( \exp\left\{\lambda \xi +(\lambda+\log(1-\lambda))\xi^2\right\}\Big) &\leq&  1.
\end{eqnarray*}
\end{lemma}
\noindent\emph{Proof.} Assume  $\xi \geq -1$ and $\lambda \in  [0, 1)$. Then $\lambda \xi \geq -\lambda> -1$.
Since the function
\begin{equation}
 f(x)=\frac{\log(1+x)-x}{ x^2/2 }, \ \ \ \ \ \ \ x > -1,
\end{equation}
is increasing in $x$, we have
\begin{eqnarray}
   \log(1+\lambda \xi )&\geq & \lambda \xi +\frac{1}{2}(\lambda \xi )^2 f(-\lambda)  \nonumber\\
   &=& \lambda \xi +\xi^2(\lambda+\log(1-\lambda)).\label{sksn}
\end{eqnarray}
Thus
\begin{eqnarray}
  \exp\left\{\lambda \xi +\xi^2(\lambda+\log(1-\lambda))\right\}
    \leq   1+ \lambda \xi .
\end{eqnarray}
Since $\mathbf{E} \xi \leq 0$, it follows that
\begin{eqnarray*}
  \mathbf{E}\Big( \exp\left\{\lambda \xi +\xi^2(\lambda+\log(1-\lambda))\right\} \Big) \leq   1,
\end{eqnarray*}
which gives the desired inequality.
\vspace{0.3cm}

\noindent\emph{Proof of Corollary \ref{th7}.}
Let $T=\min\{ k \in [1, n]: S_k\geq x \ \mbox{and} \ [S]_k\leq v^2 \}.$
Applying inequality (\ref{fhgsa2}) with $g(\lambda)=-(\lambda+\log(1-\lambda))$ and $f(\lambda)=0,$ from Lemma 
\ref{lemma17}, we obtain, for all $x,  v  >0$ and all $\lambda \in [0, 1)$,
\begin{eqnarray}
  && \mathbf{P}(  S_k\geq x \ \mbox{and} \ [S]_k\leq v^2\ \mbox{for some}\ k\in[1,n]) \nonumber\\
 &\leq&   \exp\{-\lambda x- (\lambda+\log(1-\lambda))v^2 \} . \label{fmula}
\end{eqnarray}
It is easy to see that  bound (\ref{fmula}) attains its minimum at
\begin{eqnarray}
   \lambda = \lambda(x)= \frac{x}{v^2+ x}. \label{lanbda1}
\end{eqnarray}
Substituting $\lambda=\lambda(x)$ in (\ref{fmula}), we get, for all $x,  v  >0$,
\begin{eqnarray}
&&\mathbf{P}\left( S_k  \geq x\  \mbox{and}\ [S]_k\leq v^2 \ \ \mbox{for some}\ \ k\in[1,n]\right)   \nonumber \\
&\leq& \inf_{ \lambda \in [0, 1)} \exp\{-\lambda x- (\lambda+\log(1-\lambda))v^2 \} \label{fnmsx}\\
&= & \left(1+ \frac{x}{v^2}\right)^{v^2}e^{-x}.   \label{fnmsv}
\end{eqnarray}
Using Taylor's expansion, we deduce that, for all $\lambda \in [0, 1)$,
\begin{eqnarray}
 \lambda+\log(1-\lambda) &=&-\frac{ \lambda^2}{2} \left(1 + \frac23  \lambda + \frac24  \lambda^2+...\right)  \nonumber \\
 &\geq& -\frac{ \lambda^2}{2}\left(1+\lambda+\lambda^2+...\right) \nonumber \\
 &=&  -\frac{ \lambda^2 }{2(1-\lambda)}.
\end{eqnarray}
Thus we have, for all $x,  v  >0$,
\begin{eqnarray}
   \inf_{ \lambda \in [0, 1)} \exp\{-\lambda x- (\lambda+\log(1-\lambda))v^2 \} &\leq & \inf_{ \lambda \in [0, 1)} \exp\left\{-\lambda x+ \frac{ \lambda^2v^2  }{2(1-\lambda)}\right\} \nonumber \\
    &=& B_{1}(x, 1, v) \label{copie1}\\
    &\leq& B_{2}(x, 1, v). \label{copie2}
\end{eqnarray}
Combining (\ref{fnmsv}), (\ref{copie1}) and (\ref{copie2})  together, we obtain the desired inequalities of Corollary \ref{th7}.

\vspace{0.3cm}

\noindent\emph{Proof of Corollary \ref{co2}.}
Assume $ \mathbf{E}(\xi_{i}^{l}  | \mathcal{F}_{i-1})  \leq \frac12 l!\epsilon^{l-2} \mathbf{E}(\xi_i^2 | \mathcal{F}_{i-1})$ for all $ l\geq 2$ and a constant $\epsilon \in (0, \infty)$.
Then, for all $0\leq \lambda < \epsilon^{-1}$,
\begin{eqnarray*}
 \mathbf{E}(e^{\lambda \xi _i}|\mathcal{F}_{i-1})-1 &=& \sum_{k=2}^{+\infty}\frac{\lambda^{k}}{k !} \mathbf{E}(\xi_{i}^{k} |\mathcal{F}_{i-1})  \nonumber\\
&\leq & \frac{\lambda^2}{2}\mathbf{E}(\xi_{i}^{2} |\mathcal{F}_{i-1})\sum_{k=2}^{\infty}(\lambda\epsilon)^{k-2}  \nonumber\\
 &= &\frac{\lambda^2}{2(1-\lambda\epsilon)}\mathbf{E}(\xi_{i}^{2} |\mathcal{F}_{i-1}).\nonumber
\end{eqnarray*}
Using Theorem \ref{th5}, we obtain the desired inequality (\ref{ie1g5}) with $\lambda= \overline{\lambda}$.
Since $n\log(1+\frac{t}{n}) \leq t$ for all $t\geq0,$ it follows that, for all $x,v>0$,
\begin{eqnarray}
 B_{1,n}(x, \epsilon, v)  \leq    \exp \left\{- \overline{\lambda} x +  \frac{\overline{\lambda}^2 v^2  }{2 (1-\overline{\lambda}\epsilon)} \right\}
= B_1(x, \epsilon, v).
\end{eqnarray}
This completes the proof of Corollary \ref{co2}.

\vspace{0.3cm}

\noindent\emph{Proof of Corollary \ref{co1}.}
Assume that $(\xi_{i},  \mathcal{F}_{i})_{i=1,...,n}$ are conditionally symmetric. For any $y>0$, let $\eta_i=\xi_{i}\textbf{1}_{\{|\xi_{i}|\leq y \}}$. Then $(\eta_i,  \mathcal{F}_{i})_{i=1,...,n}$ is  a sequence of bounded and conditionally symmetric martingale differences.
Using Taylor's expansion, we obtain  the following estimation of the moment generating function of $\eta_i$,
\begin{eqnarray*}
 \mathbf{E}(e^{\lambda \eta_i}|\mathcal{F}_{i-1}) &=& \mathbf{E}\left(\frac{e^{\lambda \eta_i}+e^{-\lambda \eta_i}}{2}\ \bigg|\ \mathcal{F}_{i-1} \right)  \nonumber\\
&=& 1+  \sum_{k=1}^{\infty}\frac{ \lambda^{2k}}{(2k)!} \mathbf{E}\left(  \eta_i^{2k} \ |\ \mathcal{F}_{i-1} \right).
\end{eqnarray*}
Since $|\eta_i|\leq y$, it follows that $\mathbf{E}\left(  \eta_i^{2k} \ |\ \mathcal{F}_{i-1} \right) \leq y^{2k-2}\mathbf{E}\left(  \eta_i^{2} \ |\ \mathcal{F}_{i-1} \right)$ and that
\begin{eqnarray}
 \mathbf{E}(e^{\lambda \eta_i}|\mathcal{F}_{i-1}) &\leq& 1+  \frac{\mathbf{E}\left(  \eta_i^{2} \ |\ \mathcal{F}_{i-1} \right)}{y^2} \sum_{k=1}^{\infty}\frac{ (\lambda y)^{2k}}{(2k)!} \nonumber\\
 &\leq& 1+  \frac{\mathbf{E}\left(  \eta_i^{2} \ |\ \mathcal{F}_{i-1} \right)}{y^2} \left( \cosh(\lambda y)-1\right).\label{kdagfd}
\end{eqnarray}
Set $V_k^2(y)=\sum_{i=1}^{k}\mathbf{E}(\eta_i^2 |\mathcal{F}_{i-1})$ for all $k \in [1,n]$.
Using Theorem \ref{th5}, we obtain, for all $  x, v> 0,$
\begin{eqnarray}
P_1&:=& \mathbf{P}\left(  \sum_{i=1}^k\eta_i \geq x\ \mbox{and}\ V_k^2(y) \leq v^2\ \mbox{for some}\ k \in [1,n] \right) \nonumber\\
&\leq& \inf_{  \lambda \geq 0} \exp \left\{- \lambda x + n\log \left(1+ \frac{ v^2  }{ n\,y^2}\left( \cosh(\lambda y)-1\right) \right)\right\} \label{fg25}\\
&\leq& \inf_{  \lambda \geq 0} \exp \left\{- \lambda x + \frac{ v^2  }{ y^2}\left( \cosh(\lambda y)-1\right) \right\}.\label{fg26}
\end{eqnarray}
By some simple calculations, we find that
(\ref{fg25}) and (\ref{fg26}) attain their minimums at $\underline{\lambda}$ and $\overline{\lambda}$ of Corollary \ref{co1}, respectively.
It is easy to see that
\begin{eqnarray}
&& \mathbf{P}\left(  S_k \geq x\ \mbox{and}\ V_k^2(y)\leq v^2\ \mbox{for some}\ k \in [1,n] \right)  \nonumber\\
&\leq& \mathbf{P}\left(  \sum_{i=1}^k\left(\eta_i + \xi_{i}\textbf{1}_{\{ \xi_{i} < - y \}}\right)\geq x\ \mbox{and}\ V_k^2(y) \leq v^2\ \mbox{for some}\ k \in [1,n] \right) \nonumber\\
&&+ \,\mathbf{P}\left( \sum_{i=1}^k \xi_{i}\textbf{1}_{\{ \xi_{i} > y \}}> 0\ \mbox{and}\ V_k^2(y) \leq v^2\ \mbox{for some}\ k \in [1,n] \right) \nonumber\\
&\leq& P_1\ + \, \mathbf{P}\left( \max_{1\leq i \leq n}  \xi_i > y \right) . \label{fghgh}
\end{eqnarray}
Implementing (\ref{fg25}) and (\ref{fg26}) into  (\ref{fghgh}), we get the desired inequalities (\ref{hdgh1}) and (\ref{hdgh2}).

\section{Proof of Theorem \ref{th3} and its corollaries}\label{sec5}
The proof of Theorem \ref{th3} is similar to the argument of Theorem \ref{th5}.

\noindent\emph{Proof of Theorem \ref{th3}.}
Let $T=\min\{ k \in [1, n]: S_k\geq x \ \mbox{and} \ \sum_{i=1}^{k}V_{i-1}\leq v^2 \}.$
According to (\ref{ghnda}) with $g(\lambda)\equiv 0$, we have the following
estimation,  for all $ x, v>0$,
\begin{eqnarray*}
 \mathbf{P}\bigg( S_{k} \geq x\ \mbox{and}\ \sum_{i=1}^{k}V_{i-1} \leq v^2\ \textrm{for some}\  k \in [1,n] \bigg)\ \leq  \  \sum_{k=1}^n \mathbf{E}_\lambda \Big(\exp \left\{ -\lambda x +\Psi _k(\lambda )\right\} \mathbf{1}_{\left\{ T = k \right\}}\Big),
\end{eqnarray*}
where
\begin{eqnarray*}
\Psi _k(\lambda )=  \sum_{i=1}^k \log \mathbf{E}(e^{\lambda \xi_{i}}|\mathcal{F}_{i-1}) .
\end{eqnarray*}
Using the condition $ \mathbf{E}(e^{\lambda \xi_{i}}|\mathcal{F}_{i-1})$ $\leq $$\exp \{ f(\lambda) V_{i-1}  \}$ and the fact $\sum_{i=1}^{k}V_{i-1}\leq v^2$  on the set  $\left\{ T = k \right\}$, we obtain
\begin{eqnarray*}
&& \mathbf{P}\bigg( S_{k} \geq x\ \mbox{and}\ \sum_{i=1}^{k}V_{i-1} \leq v^2\ \textrm{for some}\  k \in [1,n]  \bigg) \nonumber \\
&\leq & \sum_{k=1}^n \mathbf{E}_\lambda \Bigg(\exp \left\{ -\lambda x +f(\lambda)  \sum_{i=1}^{k}V_{i-1} \right\} \mathbf{1}_{\left\{ T = k \right\}}  \Bigg) \\
&\leq& \exp \left\{ -\lambda x +f(\lambda)  v^2 \right\} ,
\end{eqnarray*}
which gives (\ref{f1sd}) of  Theorem \ref{th3}.

In the proof  of Corollary \ref{co4}, we shall need the following two lemmas.
\begin{lemma}
\label{lemma5} If $\xi$ is a random variable satisfying $\xi \leq 1$, $\mathbf{E}\xi
\leq 0$ and $\mathbf{E}\xi^2=\sigma^2$, then, for all $\lambda > 0,$
\[
\mathbf{E}  e^{\lambda \xi}   \leq \frac{1}{1+\sigma^2} \exp\left\{-\lambda \sigma^2
\right\} + \frac{\sigma^2}{1+\sigma^2 }\exp\{\lambda\} .
\]
\end{lemma}

A proof can be found in Fan, Grama and Liu  \cite{F12}.

\begin{lemma}\label{lemma6}  Assume that $\xi$ is a random variable satisfying $\mathbf{E}\xi
\leq 0$,  $\xi \leq b$ for a constant  $b> 0$ and $\mathbf{E}\xi^2=\sigma^2$. Set
\begin{eqnarray}
s^2  = \left\{ \begin{array}{ll}
\sigma^2 , & \textrm{\ \ \ \ \ if $\sigma \geq b $}, \\
\frac{1}{4}\left(b   + \frac{\sigma^2}{b  }\right)^2, & \textrm{\ \ \ \ \ if $\sigma  < b  $}.
\end{array} \right.
\end{eqnarray}
Then, for all $\lambda > 0,$
\begin{equation}\label{ghsdfd}
\mathbf{E} e^{\lambda \xi}  \leq \exp\left\{ \frac{\lambda^2s^2}{2} \right\}.
\end{equation}
\end{lemma}
\noindent\emph{Proof.}
If $\sigma  \geq b$, by Lemma \ref{lemma5}, then, for all $t\geq0$,
\[
\mathbf{E}e^{t \xi/\sigma} \leq \frac{1}{2}\Big(e^{-t}+ e^{ t}\Big)\leq \exp\left\{ \frac{t^2}{2} \right\}.
\]
Taking $t =\lambda \sigma \geq 0$, we have
\begin{eqnarray}\label{fie1}
\mathbf{E}e^{\lambda \xi } \leq  \exp\left\{ \frac{ \lambda^2 \sigma^2}{2} \right\}=\exp\left\{\frac{\lambda^2s^2}{2} \right\}.
\end{eqnarray}
If $\sigma < b $, by Lemma \ref{lemma5}, we get, for all $t\geq0$,
\begin{eqnarray*}
\mathbf{E}e^{t \xi/b} &\leq& \frac{1}{1+\sigma^2/b^2} \exp\left\{-t \sigma^2/b^2
\right\} + \frac{\sigma^2/b^2}{1+\sigma^2/b^2 }\exp\{t\} \\
 &=& \exp\left\{ f(z)\right\},
\end{eqnarray*}
where $z=t(1+\sigma^2/b^2)$ and $f(z)=- z p + \log(1-p + p e^z)$ with $p=\frac{\sigma^2/b^2}{1+\sigma^2/b^2}$.
Since $f(0)=f'(0)=0$,
\[
f'(z)=-p+\frac{p}{p+(1-p)e^{-z}}
\]
and
\[
f''(z)=\frac{p(1-p)e^{-z}}{(p+(1-p)e^{-z})^2} \leq \frac{1}{4},
\]
we have
\[
f(z)\leq \frac{1}{8}z^2= \frac{ t^2}{8}\left(1+\frac{\sigma^2 }{ b^2}  \right)^2 \ \ \ \ \mbox{and}\ \ \ \ \ \mathbf{E}e^{t \xi/b} \leq \exp\left\{ \frac{ t^2}{8}\left(1+\frac{\sigma^2 }{ b^2}  \right)^2 \right\}.
\]
Taking $t =\lambda b\geq 0$, we obtain
\begin{eqnarray}
\mathbf{E}e^{\lambda \xi } \leq  \exp\left\{ \frac{\lambda^2}{8} \left( b +\frac{\sigma^2 }{ b }  \right)^2 \right\}=\exp\left\{\frac{\lambda^2s^2}{2} \right\}. \label{fie2}
\end{eqnarray}
Combining (\ref{fie1}) and (\ref{fie2}) together, we obtain (\ref{ghsdfd}).

\vspace{0.3cm}

\noindent\emph{Proof of Corollary \ref{co4}.}
Inequality (\ref{gnlm}) follows immediately from Lemma \ref{lemma6}. Using Theorem \ref{th3}, we obtain, for all $ x, \lambda, v>0$,
\begin{eqnarray*}
  \mathbf{P}\left(   S_k \geq x\  \mbox{and}\ \sum_{i=1}^k C_{i-1}^2 \leq v^2\ \mbox{for some}\ k \in [1, n]  \right)
  \leq   \exp\left\{- \lambda x + \frac{ \lambda^2v^2}{2}   \right\}.
\end{eqnarray*}
 Minimizing the right hand side of the last inequality with respect to $\lambda\geq 0$, we easily obtain  (\ref{f3}).

\vspace{0.3cm}

\noindent\emph{Proof of Remark \ref{eton}.}
Assume that $(\xi _i)_{i=1,...,n}$ are  independent and  symmetric. Set $$\mathcal{F}_{i}=\sigma \left\{\xi_{k}, k\leq i, \xi_j^2, 1\leq j \leq n \right\}.$$ Since $\xi _i$ is symmetric, we deduce that  $$ \mathbf{E}(\xi_i> y|\, \mathcal{F}_{i-1} ) = \mathbf{E}(\xi_i> y|\, \xi_i^2 ) =\mathbf{E}( -\xi_i > y|\, (-\xi_i)^2 )=\mathbf{E}( -\xi_i > y|\, \mathcal{F}_{i-1}  ).$$
Thus $\Big(\frac{\xi_i}{\sqrt{[S]_n}},\mathcal{F}_i\Big)_{i=1,...,n}$ are conditionally symmetric martingale differences.  For all $1\leq i\leq n$, we have
\begin{eqnarray*}
 \mathbf{E}\left(\exp\left\{ \lambda  \frac{\xi_i}{\sqrt{[S]_n}}   \right\} \Bigg|  \mathcal{F}_{i-1} \right)  = \frac12 \ \mathbf{E}\left(\exp\left\{  \lambda \frac{\xi_i}{\sqrt{[S]_n}} \right\} + \exp\left\{  -\lambda  \frac{\xi_i}{\sqrt{[S]_n}} \right\} \Bigg|  \mathcal{F}_{i-1}  \right).
\end{eqnarray*}
Using the inequality $\frac12(e^t+e^{-t}) \leq e^{t^2/2}$, we obtain, for all $\lambda\geq0$,
\begin{eqnarray*}
 \mathbf{E}\left(\exp\left\{ \lambda  \frac{\xi_i}{\sqrt{[S]_n}}  \right\} \Bigg| \mathcal{F}_{i-1} \right)
  \leq  \exp\left\{\frac{\lambda^2\xi_i^2}{2\, [S]_n } \right\}.
\end{eqnarray*}
Since $\xi_i^2$ is measurable with respect to $\mathcal{F}_{i-1},$
it follows that $$ \sum_{i=1}^{k}\mathbf{E}\left(\frac{\xi_i^2}{\,[S]_n}  \Big| \mathcal{F}_{i-1} \right)= \sum_{i=1}^{k} \frac{\xi_i^2}{\,[S]_n} \leq 1$$ for all $k \in [1, n].$
By Theorem \ref{th3} with $V_{i-1}=\frac{\xi_i^2}{\,[S]_n},$ it follows that, for all $x, \lambda\geq0$,
\begin{eqnarray}
 \mathbf{P}\left( \max_{ 1\leq k \leq n  } \frac{S_k}{\sqrt{[S]_n} } \geq x \right)
  \leq  \exp \left\{ -\lambda x +\frac{\lambda^2}{2} \right\}. \label{fjna}
\end{eqnarray}
The right hand side of the last inequality attends its minimum at $\lambda=x$. Substituting $\lambda=x$ into (\ref{fjna}),
we easily get (\ref{fdgfgfdgh}) of Remark \ref{eton}.

\section{Proof of Theorems \ref{thlin} - \ref{th6}}\label{endsec}
We make use of Corollary \ref{co2} to prove Theorem \ref{thlin}.\\
\emph{Proof of Theorem \ref{thlin}.} From (\ref{ine29}) and (\ref{ine30}), it is easy to see that
\begin{equation}\label{sdgvf1}
\theta_n -\theta = \sum_{k=1}^n\frac{ \phi_{k} \varepsilon_k}{\sum_{k=1}^n \phi_{k}^2}. \nonumber
\end{equation}
For any $i=1,...,n$, set
\begin{eqnarray}\label{sdgvf2}
\xi_i= \frac{ \phi_{i} \varepsilon_i}{ \sigma \sqrt{\sum_{k=1}^n \phi_{k}^2}}\ \ \ \  \textrm{and} \ \ \ \ \mathcal{F}_{i} = \sigma \Big( \phi_{k}, \varepsilon_k, 1\leq k\leq i,\  \phi_{k}^2, 1\leq k\leq n \Big).
\end{eqnarray}
Then $(\xi _i,\mathcal{F}_i)_{i=1,...,n}$ is a sequence of martingale differences and satisfies
$$ \frac{(\theta_n -\theta)\sqrt{\sum_{k=1}^n \phi_{k}^2} } { \sigma}  =\sum_{i=1}^n\xi_i.$$
Notice that $$\langle S\rangle_n = \sum_{i=1}^{n} \frac{ \phi_{i}^2 }{ \sigma^2 (\sum_{k=1}^n \phi_{k}^2)} \mathbf{E}(\varepsilon_i^2 | \mathcal{F}_{i-1} ) =\sum_{i=1}^{n} \frac{ \phi_{i}^2 }{  \sum_{k=1}^n \phi_{k}^2 }  = 1, $$
and that
\begin{eqnarray} \label{fgdsasas}
|\mathbf{E}(\xi_i^k | \mathcal{F}_{i-1} )| &=&  \frac{ \phi_{i}^2 }{ \sigma^k (\sum_{k=1}^n \phi_{k}^2)} \mathbf{E}\bigg( \Big(\frac{ \phi_{i}  }{  \sqrt{\sum_{k=1}^n \phi_{k}^2}}\Big)^{k-2} \varepsilon_i^k  \bigg| \mathcal{F}_{i-1} \bigg) \nonumber\\
&\leq&  \frac{ \phi_{i}^2 \epsilon_1^{k-2} }{ \sigma^k (\sum_{k=1}^n \phi_{k}^2)} \mathbf{E} \Big(   \varepsilon_i^k  \Big| \mathcal{F}_{i-1}  \Big) \nonumber\\
&\leq& \frac12 k! \, \epsilon_2^{k-2} \frac{ \phi_{i}^2 \epsilon_1^{k-2} }{ \sigma^{k-2}(\sum_{k=1}^n \phi_{k}^2)} \nonumber\\
&=&  \frac12 k! \, \epsilon ^{k-2} \mathbf{E}(\xi_i^2 | \mathcal{F}_{i-1} ). \nonumber
 \end{eqnarray}
Applying Corollary \ref{co2}  to $(\xi _i,\mathcal{F}_i)_{i=1,...,n}$, we obtain the claim of Theorem \ref{thlin}.

\vspace{0.3cm}

\noindent \emph{Proof of Theorem \ref{dssaf}.}     It is easy to see that the martingale differences
 $(\xi _i,\mathcal{F}_i)_{i=1,...,n}$, defined by (\ref{sdgvf2}),  satisfy
\begin{eqnarray}
 \xi_i  \leq  U_{i-1}:= \frac{ |\phi_{i}| \epsilon }{ \sigma \sqrt{\sum_{k=1}^n \phi_{k}^2}}\ \ \ \ \ \textrm{and} \ \ \ \  \mathbf{E}(\xi_i^2|\mathcal{F}_{i-1})=\frac{ \phi_{i}^2 }{  \sum_{k=1}^n \phi_{k}^2 }.
\end{eqnarray}
Applying Corollary \ref{co4}  to $(\xi _i,\mathcal{F}_i)_{i=1,...,n}$, we obtain the desired inequality.

\vspace{0.3cm}

\noindent\emph{Proof of Theorem \ref{thds}.} From (\ref{ine29}) and (\ref{ine30}), it is easy to see that
\begin{equation}\label{sdgvf1}
(\theta_n -\theta )\sum_{k=1}^n \phi_{k}^2 = \sum_{k=1}^n\phi_{k} \varepsilon_k  . \nonumber
\end{equation}
For any $i=1,...,n$, set
\begin{eqnarray}
\xi_i=\phi_{i} \varepsilon_i   \ \ \ \textrm{and} \ \ \ \ \  \mathcal{F}_{i} = \sigma \Big( \phi_{k}, \varepsilon_k, 1\leq k\leq i, \phi_{i+1}  \Big).
\end{eqnarray}
Then $(\xi _i,\mathcal{F}_i)_{i=1,...,n}$ is a sequence of martingale differences  and satisfies
$$ \mathbf{E}(e^{ \lambda \xi_{i} }|\mathcal{F}_{i-1})   \leq e^{c  |\lambda \phi_{i}|^\alpha  }\ \textrm{ for all }i \in [1, n] .$$
Applying Theorem  \ref{th3}   to $(\xi _i,\mathcal{F}_i)_{i=1,...,n}$, we obtain, for all $x, \lambda,  v \geq 0$,
\begin{eqnarray}\label{scsasd}
 \mathbf{P}\Big(   \pm (\theta_n -\theta)  \sum_{k=1}^n \phi_{k}^2  \geq x  \  \textrm{and} \   \sum_{k=1}^n   |\phi_{k}|^\alpha \leq v^\alpha  \Big)  &\leq&  \exp\Big\{ - \lambda x + c \lambda^\alpha v^\alpha   \Big\}.
\end{eqnarray}
The right hand side of the last inequality takes its minimum at $$\lambda= \lambda(x)=\Big(\frac{x}{c \, \alpha \, v^\alpha } \Big)^{\frac{ 1}{\alpha -1}} .$$
Substituting $\lambda= \lambda(x)$ into (\ref{scsasd}), we obtain the desired inequality.

\vspace{0.3cm}

\noindent\emph{Proof of Theorem \ref{th6}.} From (\ref{scfso1}) and (\ref{scfso2}), it is easy to see that
\begin{equation}\label{sdgvf1}
(\theta_n' -\theta) \sum_{k=1}^n X_{k-1}^2= \sum_{k=1}^n X_{k-1} \varepsilon_k . \nonumber
\end{equation}
For any $i=1,...,n$, set
\begin{eqnarray}
\xi_i=X_{i-1} \varepsilon_i   \ \ \ \textrm{and} \ \ \ \ \  \mathcal{F}_{i} = \sigma \Big( X_{0}, \varepsilon_k, 1\leq k\leq i \Big).
\end{eqnarray}
Then $(\xi _i,\mathcal{F}_i)_{i=1,...,n}$ is a sequence of martingale differences and satisfies
$$ |\xi_i|   \leq U_{i-1}:= X_{i-1}\epsilon \ \ \ \ \ \textrm{and} \ \ \ \ \    \mathbf{E} (( X_{k-1} \varepsilon_{i})^2|\mathcal{F}_{i-1} ) =  X_{k-1}^2  \mathbf{E}  \varepsilon_{i}^2     \leq  X_{k-1}^2 \sigma^2  .$$
Applying Corollary \ref{co4} to $(\xi_i, \mathcal{F}_{i})$, we obtain the desired inequality.

\ACKNO{We would like to thank M. Ledoux, D. Chafa\"{\i} and the referee for their helpful comments and suggestions.
Fan was partially supported by the Post-graduate Study Abroad Program sponsored by China Scholarship Council.
Liu was partially supported by the  National Natural Science Foundation of China (no. 11171044, no. 11101039, and no. 11401590).
 }

\end{document}